\input amstex
\input amsppt.sty
\magnification=\magstep1  \hsize=30truecc  \vsize=22.2truecm
\baselineskip=16truept  \parskip=2pt \NoBlackBoxes  \TagsOnRight
\pageno=1 \nologo
\def\Z{\Bbb Z}          
\def\l{\left}  \def\r{\right}  \def\bg{\bigg}  \def\({\bg(}  \def\[{\bg\lfloor}  \def\){\bg)}  \def\]{\bg\rfloor}
\def\t{\text}  \def\f{\frac}  
  
\def\eq{\equiv}  \def\={\;=\;} \def\weq{\;\eq\;} \def\+{\,+\,}

\def\bi{\binom}   \def\cs{\cdots}
\def\ls{\leqslant}  
\def\mo{\roman{mod}}

\def\Ack{\medskip\noindent {\bf Acknowledgments}}
\topmatter
\title On a curious property of Bell numbers\endtitle
\author Zhi-Wei Sun and Don Zagier\endauthor
\leftheadtext{Zhi-Wei Sun and Don Zagier} \rightheadtext{On a
curious property of Bell numbers}
\abstract In this paper we derive congruences  expressing Bell numbers and derangement
numbers in terms of each other modulo any prime.
\endabstract
\thanks 2010 {\it Mathematics Subject Classification}.\,Primary 11B75;
Secondary 05A15, 05A18, 11A07.
\newline\indent The first author is supported by the National Naturaal Science
Foundation (grant 10871087) and the Overseas Cooperation Fund (grant
10928101) of China.
\endthanks
\endtopmatter
\document

\def\AA{1}  \def\AB{2}  \def\AC{3}  
 \def\BX{4}  \def\BA{5} \def\BB{6}       

\heading{1. Introduction}\endheading

Let $B_n$ denote the $n$th Bell number, defined as the number of
partitions of a set of cardinality $n$ (with $B_0=1$). In 1933
Touchard [T] proved that for any prime $p$ we have
$$ B_{p+n} \weq B_n+B_{n+1}\ (\mo\ p)\ \ \ \t{for all}\ n=0,1,2,\ldots.  \tag \AA$$
Thus it is natural to look at the numbers $B_n\pmod p$ for $n<p$.
In~[S], the first author discovered experimentally that for a fixed
positive integer~$m$ the sum $\sum_{n=0}^{p-1}B_n/(-m)^n$ modulo a
prime $p$ not dividing $m$ is independent of the prime $p$, a
typical case being
$$  \sum_{n=0}^{p-1}\frac{B_n}{(-8)^n} \weq -1853 \pmod p \qquad\t{for all primes $p\ne2\,$.} $$
In this note we will prove this fact and give some related results.

Our theorem involves another combinatorial quantity, the derangement
number $D_n$, defined either as the number of fixed-point-free
permutations of a set of cardinality~$n$ (with $D_0=1$) or by the
explicit formula
$$\f{D_n}{n!}=\sum_{k=0}^n\f{(-1)^k}{k!} \qquad(n=0,\,1,\,2,\,\dots)\;.  \tag\AB$$

\proclaim{Theorem 1} For every positive integer $m$ and any prime
$p$ not dividing $m$ we have
$$\sum_{0<k<p}\f{B_k}{(-m)^k}\;\eq\; (-1)^{m-1}D_{m-1}\pmod p\;.   \tag\AC $$
\endproclaim

Using $\sum_{0<m<p}(-m)^{n-k}\eq-\delta_{n,k}\pmod p$ for
$k,\,n\in\{1,\dots,p-1\}$, we immediately obtain a dual formula for
$B_n$ $(n<p)$ in terms of ~$D_0,\ldots,D_{p-2}$.

\proclaim{Corollary}
Let $p$ be any prime. Then for all $n=1,\ldots, p-1$ we have
$$ B_n\eq\sum_{m=1}^{p-1}(-1)^mD_{m-1}\,(-m)^n \pmod p \;. $$
\endproclaim

For the reader's convenience we give a small table of values of
$B_n$~and~$D_n$.
$$\vbox{\offinterlineskip\halign{$\hfil#\hfil$&\ \vrule#\ &$\hfil#\ \hfil$&$\hfil\ #\ \hfil$
  &$\hfil\ #\ $&$\hfil\ #\ $&$\hfil\ #\ $&$\hfil\ #\ $&$\hfil\ #\ $&$\hfil\ #\ $&$\hfil\ #\ $&$\hfil\ #$\cr
n&&&0&1&2&3&4&5&6&7&8 \cr &height3pt& &&&&&&&&\cr
\noalign{\hrule}\omit &height2pt&\omit&&&&&\cr
B_n&&&1&1&2&5&15&52&203&877&4140 \cr \omit
&height3pt&\omit&&&&&&&\cr     \noalign{\hrule}\omit
&height2pt&\omit&&&&&\cr D_n&&&1&0&1&2&9&44&265&1854&14833\cr \omit
&height4pt&\omit&&&&&&&&\cr  \omit &height3pt&\omit&&&&&&&\cr}}$$

We will show Theorem 1 in the next section and derive an extension
of Theorem 1 in Section 3.

\heading{2. Proof of Theorem 1}\endheading We first observe that it
suffices to prove (\AC) for $0<m<p$, since both sides are periodic
in $m\pmod p$ with period~$p$. For the left-hand side this is
obvious and for the right-hand side it follows from~(\AB), which
gives the expression $(-1)^nD_n\=\sum_{r=0}^\infty(-1)^r
n(n-1)\cs(n-r+1)$ for $D_n\pmod p$ as a terminating infinite series
of polynomials in~$n$.

We will prove (\AC) for $0<m<p$ by induction on~$m$. Denote by~$S_m$
the sum on the left-hand side of~(\AC), where we consider the
prime~$p$ as fixed and omit it from the notation.  Since
$D_n=nD_{n-1}+(-1)^n$ for $n=1,2,3,\ldots$ (obvious from~(\AB)), we
have to prove the two formulas
$$  S_1\weq 1\!\pmod p\,, \qquad m\,S_m\weq S_1\,-\,S_{m+1}\!\pmod p\,. \tag\BX$$

Recall that the Bell numbers can be given by the generating function
$$\sum_{n=0}^\infty B_n\f{x^n}{n!}\=e^{e^x-1} \tag\BA$$
equivalent to the well-known closed formula
$$B_n=\dfrac1e\,\sum\limits_{r=0}^\infty\dfrac{r^n}{r!}\,.$$
Since the function $y=e^{e^x-1}$ satisfies $y'=e^xy$, this also
gives the recursive definition
$$ B_0\=1\,, \qquad B_{n+1}\=\sum_{k=0}^n\bi nk\, B_k\qquad \t{for all $n\ge0\,$.}  \tag\BB $$
This recursion is the key ingredient in proving (\BX).

For the first formula in (\BX) we use (\BB) with $n=p-1$ to obtain
$$S_1\=\sum_{k=1}^{p-1}(-1)^k\,B_k \weq \sum_{k=1}^{p-1}\bi{p-1}k\,B_k \=B_p-B_0\pmod p\,,$$
so it suffices to prove that $ B_p\eq2\pmod p$.  This is a special
case of Touchard's congruence~(\AA), but can also be seen by writing
(\BA) in the form
$$\sum_{n=0}^\infty B_n\f{x^n}{n!}\=e^x \+ \sum_{1<r<p}\f{\bigl(e^x-1\bigr)^r}{r!}
  \+ \f{x^p}{p!}\+\t O\bigl(x^{p+1}\bigr)$$
to get $B_p/p!\,=\,1/p!\+\t{($p$-integral)}\+1/p!\,$.

Now using Fermat's little theorem we have
$$ \align -m\,S_m &\weq \sum_{n=0}^{p-2} (-m)^{p-1-n}\sum_{k=0}^n \bi nk\,B_k
 \\ &\weq \sum_{k=0}^{p-2}(-1)^k\,B_k\,\sum_{r=0}^{p-k-2} \bi{p-k-1}r\,m^{p-k-1-r} \qquad\quad(r=n-k) \\
 & \weq \sum_{k=0}^{p-2} (-1)^kB_k\,\bigl((m+1)^{p-1-k}-1\bigr) \weq S_{m+1}\,-\,S_1 \pmod p  \endalign$$
for $1\le m\le p-2$.  This completes the proof of (\BX) and the theorem.

\heading{3. An extension of Theorem 1}\endheading

 Recall that for nonnegative integers $n$ and $k$ the Stirling number $S(n,k)$ of the second kind
 is the number of ways to partition a set of $n$ elements into $k$ groups.
 Obviously $$B_n\=\sum_{k=0}^n S(n,k).$$
 The Touchard polynomial $T_n(x)$ of degree $n$ is given by
 $$T_n(x)\=\sum_{k=0}^nS(n,k)x^k.\tag7$$ Note that $T_n(1)\=B_n$. Similar to the recursion for Bell numbers, we have the recursion
 $$T_{n+1}(x)\=x\sum_{k=0}^n\bi nkT_k(x).\tag8$$

 Let $p$ be a prime and let $\Z_p$ be the ring of $p$-adic integers. For two polynomials $P(x),Q(x)\in\Z_p[x]$,
 by $P(x)\weq Q(x)\ (\mo\ p)$
 we mean that the corresponding coefficients of $P(x)$ and $Q(x)$ are congruent modulo $p$.

 Our next theorem is a further generalization of Theorem 1.

 \proclaim{Theorem 2} For every positive integer $m$, we have
 $$(-x)^m\sum_{0<n<p}\f{T_n(x)}{(-m)^n}\weq-x^p\sum_{l=0}^{m-1}\f{(m-1)!}{l!}(-x)^l\pmod p\tag9$$
 for any prime $p$ not dividing $m$.
 \endproclaim

 As a consequence, if $x$ is a $p$-adic integer not divisible by $p$, then
  $$\sum_{0<n<p}\f{T_n(x)}{(-m)^n}\weq\f1{(-x)^{m-1}}\sum_{k=0}^{m-1}\f{(m-1)!}{k!}(-x)^k\pmod p.\tag10$$
  In particular,
 $$\align\sum_{0<n<p}\f{T_n(x)}{(-2)^n} &\weq\f {x-1}x\pmod p \qquad\t{for $p\ne2$}, \\
 \sum_{0<n<p}\f{T_n(x)}{(-3)^n}&\weq\f {x^2-2x+2}{x^2}\pmod p\qquad\t{for $p\ne3$}\,,\\
 \sum_{0<n<p}\f{T_n(x)}{(-4)^n}&\weq\f {x^3-3x^2+6x-6}{x^3}\pmod p \qquad\t{for $p\ne2$}\,.
 \endalign$$

\medskip

Although we can show Theorem 2 via a slight modification of the proof of Theorem 1, below we prove Theorem 2 by a new approach.

\medskip

\noindent{\it Proof of Theorem 2}. Observe that
$$\sum_{n=1}^{p-1}\f{T_n(x)}{(-m)^n}\=\sum_{n=1}^{p-1}\f{\sum_{k=1}^nS(n,k)x^k}{(-m)^n}\=\sum_{k=1}^{p-1}x^k\sum_{n=1}^{p-1}\f{S(n,k)}{(-m)^n}.$$
It is known that
$$S(n,k)\=\f1{k!}\sum_{j=0}^k\bi kj(-1)^{k-j}j^n\quad\t{for all}\ n,k=0,1,2,\ldots.$$
Thus
$$\sum_{n=1}^{p-1}\f{T_n(x)}{(-m)^n}\=\sum_{k=1}^{p-1}\f{x^k}{k!}\sum_{n=1}^{p-1}\sum_{j=1}^k\bi kj(-1)^{k-j}\l(-\f jm\r)^n.$$
For each $j\in\{1,\ldots,p-1\}$, if $p\mid m+j$ then
$$\sum_{n=1}^{p-1}\l(-\f jm\r)^n\weq\sum_{n=1}^{p-1}1\weq-1\pmod p\,,$$
if $p\nmid m+j$ then
$$\sum_{n=1}^{p-1}\l(-\f jm\r)^n\=\sum_{n=0}^{p-1}\l(-\f jm\r)^n-1\=\f{(-j/m)^p-1}{-j/m-1}-1\weq0\pmod p$$
with the help of Fermat's little theorem.

 Let $r$ denote the least positive residue of $-m$ modulo $p$. By the above,
$$\align\sum_{n=1}^{p-1}\f{T_n(x)}{(-m)^n}\=&\sum_{k=1}^{p-1}\f{x^k}{k!}\sum_{j=1}^{p-1}\bi kj(-1)^{k-j}\sum_{n=1}^{p-1}\l(-\f jm\r)^n
\\\weq&\sum_{k=1}^{p-1}\f{x^k}{k!}\bi kr(-1)^{k-r}(-1)=\f{(-1)^{r+1}}{r!}\sum_{k=r}^{p-1}\f{(-x)^k}{(k-r)!}\pmod p.
\endalign$$
Therefore
$$(-x)^m\sum_{n=1}^{p-1}\f{T_n(x)}{(-m)^n}\weq\f{(-1)^{r+1}}{r!}\sum_{l=m+r-p}^{m-1}\f{(-x)^{p+l}}{(p+l-m-r)!}\pmod p.$$
So it remains to show that
$$\f{(m-1)!}{l!}\weq\cases (-1)^{r+1}/(r!(p+l-m-r)!)\pmod p&\t{if}\ m+r-p\ls l<m,
\\0\pmod p&\t{if}\ 0\ls l<m+r-p.\endcases$$

If $0\ls l<m+r-p$, then we have $(m-1)!/l!\weq0\ (\mo\ p)$ since
$l<m+r-p<m$ and $m+r\weq0\pmod p$.

Now suppose that $m+r-p\ls l<m$. Then
$$\align&\f{(-1)^{r+1}}{(p+l-m)!}\bi{p+l-m}r
\\\weq& \f{(-1)^{r+1}\prod_{s=1}^{m-l-1}(p-s)}{(p-1)!}\bi{l-m}r
\\\weq&(-1)^{m-l-1}(m-l-1)!\bi{m-l+r-1}{m-l-1}\ (\t{by Wilson's theorem})
\\\weq&(-1)^{m-l-1}(m-l-1)!\bi{-l-1}{m-l-1}=\f{(m-1)!}{l!}\pmod p\ \ (\t{as}\ p\mid m+r).
\endalign$$

In view of the above, we have completed the proof of Theorem 2. \qed

\Ack. The joint work was done during the authors' visit to the
National Center for Theoretical Sciences (Hsinchu, Taiwan) during
August 1--8, 2010. Both authors are indebted to Prof. Winnie
Wen-Ching Li for the kind invitation and the center for the
financial support.

\address Department of Mathematics, Nanjing University,
   Nanjing 210093, People's Republic of China\endaddress
 \email {zwsun\@nju.edu.cn}\ {\it Homepage}:\ {\tt http://math.nju.edu.cn/$\sim$zwsun}
 \endemail

\address Max-Planck-Institut f\"ur Mathematik, Bonn 53111, Germany\endaddress
\email {don.zagier\@mpim-bonn.mpg.de}\endemail

\enddocument